\documentclass[12pt]{amsart}
\usepackage{amsfonts}
\usepackage{amssymb}
\theoremstyle{definition}

\theoremstyle{remark}

\newcommand{\ds}{\displaystyle}

\begin{document}

\centerline{\bf ANNUAIRE \ DE \ L'UNIVERSIT\'E \ DE \ SOFIA \ "KLIMENT \ OHRIDSKI"}

\vspace{0.2in}
\centerline{\bf FACULT\'E DE MATH\'EMATIQUES ET M\'ECANIQUE}

\vspace{0.2in}
\centerline{\it Tome 78  \qquad\qquad\qquad\qquad\qquad Livre 1 -- Math\'ematiques   \qquad\qquad\qquad\qquad\qquad   1984}

\vspace{0.6in}

\vspace{0.2in}

\centerline{\large\bf CURVATURE PROPERTIES AND ISOTROPIC  }
\centerline{\large\bf PLANES OF RIEMANNIAN AND ALMOST HERMITIAN}
\centerline{\large\bf MANIFOLDS OF INDEFINITE METRICS}

\vspace{0.3in}
\centerline{\bf Adrijan Borisov, Georgi Ganchev, Ognian Kassabov}

\vspace{0.3in}
{\sl  The behaviour of the curvature tensor on isotropic planes on a Riemannian manifold
of indefinite metric is closely related to the curvature properties of the manifold.
The importance of isotropic planes in questions concerning curvature properties of a
Riemannian manifold $M$ with indefinite metric firstly has been shown by Dajczer and
Nomizu in [1]. Namely, they have shown that $M$ is of constant sectional curvature 
if and only if the curvature tensor vanishes on all the isotropic planes. 
In this paper we study two types of isotropic planes: weakly isotropic and strongly
isotropic planes. We prove that a Riemannian manifold of indefinite metric is conformally 
flat if and only if its curvature tensor vanishes on all the strongly isotropic
planes. We specialize the plane axiom for Riemannian manifolds of indefinite metrics.
We show that manifolds satisfying plane axiom of weakly (strongly) isotropic planes
are of constant sectional curvature (conformaly flat). Further we study analogous
problems on almost Hermitian manifolds of indefinite metrics taking into account both
structures: the metric and the almost complex structure.}

\vspace{0.5in}
\centerline{\bf 1. RIEMANNIAN MANIFOLDS OF INDEFINITE METRICS }

\vspace{0.2in}
{\bf Preliminaries.} Let $M$ be a Riemannian manifold with indefinite metric $g$
of signature $(s,n-s)$, i.e. the tangent space $T_pM$ in $p\in M$ is isometric to 
$ {\bf R}_s^n$ with the inner product:
$$
	<x,y>=-\sum_{i=1}^s x^iy^i+\sum_{j=s+1}^n x^jy^l\ .
$$
The curvature operator on  $M$  is given by
$$
	R(X,Y)=[\nabla_X,\nabla_Y]-\nabla_{[X,Y]} \ ,
$$
for arbitrary vector fields $X,\, Y$ on $M$. The Ricci tensor and the scalar curvature 
of $M$ are denoted by $\rho$ and $\tau$ respectively. A 2-plane (2-dimensional subspace)
$\alpha$ of $T_pM$ is said to be nondegenerate, weakly isotropic or strongly isotropic
if the restriction of $g$ on $\alpha$ is of rank 2, 1 or 0 respectively. The sectional 
curvature of a nondegenerate 2-plane $\alpha$ in $T_pM$ with a basis $\{x,y\}$ is given by
$$
	K(\alpha)=K(x,y)=\frac{R(x,y,y,x)}{g(x,x)g(y,y)-g^2(x,y)}\ .
$$
A pair $\{ x,a \}$ of tangent vectors at a point $p\in M$ is said to be orthonormal of
signature $(+,-)$ if $g(x,x)=1$, $g(a,a)=-1$, $g(x,a)=0$. Orthonormal pairs of signature 
$(+,+)$ or $(-,-)$ are determined analogously. Orthonormal quadruples of certain 
signature are determined in a similar way.

Firstly, the particular importance of degenerate planes in pseudo-Riemannian geometry
has been shown by the theorem:

\vspace{0.07cm}
Theorem A [1]. If $R(x,y,y,x)=0$, whenever span$\{x,y\}$ is weakly isotropic, then all the 
nondegenerate 2-planes have the same sectional curvature, i.e. the manifold is of constant 
sectional curvature.

\vspace{0.07cm}
In fact, this theorem is formulated in [1] for arbitrary degenerate 2-planes, but the proof
is also valid for the proposition in the above given form. 

A tensor $T$ of type (0,4) over $T_pM$, $p\in M$, is said to be curvature-like tensor if it 
has the following properties:

a) $T(x,y,z,u)=-T(y,x,z,u)$,

b) $T(x,y,z,u)=-T(x,y,u,z)$,

c) $T(x,y,z,u)+T(y,z,x,u)+T(z,x,y,u)=0$\\
for all $x,\,y,\,z,\,u$ in $T_pM$.

We shall use Theorem 1, a) from [1] in the following form:

\vspace{0.07cm}
Theorem B. Let $T$ be a curvature-like tensor over $T_pM$, $p\in M$. If
$T(x,y,z,x)=0$, whenever $\{x\,y\}$ is an orthonormal pair of signature
$(+,-)$ and $g(x,z)=g(y,z)=0$, then all the nondegenerate 2-planes have the 
same sectional curvature with respect to $T$.

\vspace{0.2cm}
{\bf Conformally flat Riemannian manifolds and isotropic planes.} The conformal 
curvature tensor $C$ on a Riemannian manifold is given by
$$
	C=R-\frac{1}{n-2}\varphi+\frac{\tau}{(n-1)(n-2)}\pi_1 \ ,
$$
where
$$
	\varphi(x,y,z,u)=g(y,z)\rho(x,u)-g(x,z)\rho(y,u)+g(x,u)\rho(y,z)-g(y,u)\rho(x,z)\ ,
$$
$$
	\pi_1(x,y,z,u)=g(y,z)g(x,u)-g(x,z)g(y,u)\ .
$$
It is well known, that, if $n>3$, then $M$ is conformally flat if and only if $C$
vanishes.

The main result in this section is the following theorem, which clarifies
the relation between the conformal curvature tensor and the behaviour of $R$
on strongly isotropic 2-planes.

\vspace{0.07cm}
T\,h\,e\,o\,r\,e\,m \, 1. Let $s\ge 2$, $n-s\ge 2$. Then $M$ is conformally flat if 
and only if $R(\xi,\eta,\eta,\xi)=0$, whenever span$\{\xi,\,\eta\}$ is strongly
isotropic.

{\it Proof}. Let $\{x,\,y,\,a,\,b\}$ be an arbitrary orthonormal quadruple of signature
$(+,+,-,-)$. Then
$$
	R(x+a,y+b,y+b,x+a)=0
$$  
implies
$$
	R(x,y,b,a)+R(x,b,y,a)=0 \ .
$$
Hence, using the first Bianchi identity, we find
$$
	R(x,y,b,a)=0   \ .\leqno (1)
$$
Replacing $\{x,\,b\}$ by $\{(x+tb)/\sqrt{1-t^2},\,(tx+b)/\sqrt{1-t^2})\}$, $|t|<1$, in (1),
we obtain
$$
	R(x,y,a,x)+R(b,y,a,b)=0 \ .
$$
From here it is not difficult to derive
$$
	(n-2)R(x,y,a,x)-\rho(y,a)=0  \ .   \leqno (2)
$$
Replacing $y$ by $(y+tb)/\sqrt{1-t^2}$, $|t|<1$, in (2), we find
$$
	(n-2)R(x,b,a,x)-\rho(b,a)=0 \ .  \leqno (3)
$$
From (2), (3) and Theorem B it follows $C=0$.

The inverse is a simple verification.

\vspace{0.07cm}
Using Theorem 1 and the continuity of $R$, in a similar way we have

\vspace{0.1cm}
T\,h\,e\,o\,r\,e\,m \, 2. Let $M$ ($s\ge 2$, $n-s\ge 2$) be a Riemannian manifold of 
indefinite metric. Then the following conditions are equivalent:

1) $M$ is conformally flat;

2) $R(x,y,a,b)=0$, whenever  $\{x,\,y,\,a,\,b\}$ is an orthonormal quadruple of signature
$(+,+,-,-)$;

3) $K(x,y)+K(a,b)=K(x,a)+K(y,b)$, whenever  $\{x,\,y,\,a,\,b\}$ 
is an orthonormal quadruple of signature $(+,+,-,-)$.
 
\vspace{0.07cm}
This proposition specifies the caracterizations of conformally flat manifolds
given by Schouten [2, p.307] and Kulkarni [3]. 

\vspace{0.2cm}
{\bf Plane axioms in the Riemannian geometry of indefinite metric.} Let $N$
(dim $N=m\ge 3$) be a differentiable manifold with a linear connection of zero torsion. 
In general, $N$ is said to satisfy the axiom of $r$-planes $(2\le r<m)$, if, for each 
point $p$ and for any $r$-dimensional subspace $E$ of $T_pM$, there exists a totally 
geodesic submanifold $N'$ containing $p$ such that $T_pN'=E$. 

Now we specialize the general plane axiom to the considered manifolds.

A Riemannian manifold $M$ of indefinite metric is said to satisfy the axiom of 
weakly (strongly) isotropic 2-planes, if, for each point $p$ and for any weakly 
(strongly) isotropic 2-plane $E$ in $T_pM$, there exists a totally geodesic submanifold 
$M'$ containing $p$ such that $T_pM'=E$ and all the tangent spaces of $M'$ are
weakly (strongly) isotropic in $M$.

\vspace{0.07cm}
T\,h\,e\,o\,r\,e\,m \, 3. Let $M$ ($n\ge 3$) be a Riemannian manifold of indefinite 
metric. If $M$ satisfies the axiom of weakly isotropic 2-planes, then $M$ is of 
constant sectional curvature.

{\it Proof.} Let $\xi$ be an isotropic vector in $T_pM$ and $x\perp \xi$, such that
$g(x,x)=1$ or $g(x,x)=-1$. By the conditions of the theorem, there exists a totally 
geodesic submanifold $M'$ of $M$ through $p$ such that $T_pM'=$span$\{\xi,\,x\}$.
Since $M'$  is totally geodesic and the connection is symmetric, $R(\xi,x)x$ is
tangent to $M'$ at $p$. This implies $R(\xi,x,x,\xi)=0$. According to Theorem A, $M$
is of constant sectional curvature.

\vspace{0.07cm}
T\,h\,e\,o\,r\,e\,m \, 4. Let $M$ ($s\ge2,\,n-s\ge2$) be a Riemannian manifold of 
indefinite metric. If $M$ satisfies the axiom of strongly isotropic 2-planes, then 
$M$ is conformaly flat.

{\it Proof.} Let $\xi,\,\eta$ be arbitrary isotropic perpendicular vectors in
$T_pM$, $p\in M$. By the condition of the theorem, there exists a totally geodesic 
submanifold  $M'$ of $M$ containing $p$ such that $T_pM'=$span$\{\xi,\,\eta\}$. Then
$R(\xi,\eta)\eta$ is in $T_pM'$. Hence, $R(\xi,\eta,\eta,\eta,\xi)=0$. Applying
Theorem 1, we obtain $M$ is conformally flat.

\vspace{0.5cm}

\centerline{\bf 2. ALMOST HERMITIAN MANIFOLDS OF INDEFINITE METRICS}

\vspace{0.5cm}
{\bf Preliminaries.} In this part $M$ will stand for an almost Hermitian manifold
with indefinite metric $g$ of signature $(2s,2(n-s))$ and almost complex structure 
$J$, i.e.
$$
	J^2=-id\ , \quad g(JX,JY)=g(X,Y) 
$$ 
for arbitrary vector fields $X,\,Y$ on $M$.

If $\nabla J=0$, then $M$ is a Kaehlerian manifold of indefinite metric.

A 2-plane $\alpha$ in $T_pM$ is said to be holomorphic (antiholomorphic) if
$J\alpha=\alpha$ ($J\alpha \perp \alpha, \ J\alpha \ne \alpha$). A pair $\{x,\,y\}$
of vectors in $T_pM$ is said to be holomorphic (antoholomorphic) if span$\{x,\,y\}$
is a holomorphic (antiholomorphic) 2-plane. Antiholomorphic triples and
quadruples are determined in a similar way.

The manifold $M$ is said to be of pointwise constant holomorphic (antiholomorphic) 
sectional curvature if in each point $p\in M$ all the nondegenerate holomorphic
(antiholo\-morphic) 2-planes have the same sectional curvature, which is a function 
of the point. As in the definite case, if $M$ is a Kaehler manifold of pointwise 
constant holomorphic (antiholomorphic) sectional curvature, then $M$ is of constant 
holomorphic sectional cur\-vature.

Let $\{ e_i,\ i=1,\hdots,2n\}$, be an orthonormal basis of $T_pM$. The Ricci $*$-tensor
and the scalar $*$-curvature are given by
$$
	\rho^*(y,z)=\sum_{i=1}^{2n}g(e_i,e_i)R(e_i,y,Jz,Je_i) \ ,
$$
$$
	\tau^*=\sum_{i=1}^{2n}g(e_i,e_i)\rho^*(e_i,e_i) \
$$
respectively.

The proof of the next proposition, given in [4] is also valid in the indefinite case.

\vspace{0.1cm}
Theorem  C. Let $T$ be a curvature-like tensor of type (0,4) over $T_pM$. If

1) $T(x,Jx,Jx,x)=0$, for an arbitrary vector $x$ in $T_pM$,

2) $T(x,y,y,x)=0$, whenever span$\{x,\,y\}$ is an antiholomorphic 2-plane,

3) $T(x,Jx,y,x)=0$, whenever span$\{x,\,y\}$ is an antiholomorphic 2-plane, \\
then $T=0$.

\vspace{0.2cm}
{\bf Almost Hermitian manidolds of pointwise constant antiholomorphic sectional curvature.}
In this section we give an analogue of Theorem A for almost Hermitian manifolds of
indefinite metric.

\vspace{0.07cm}
T\,h\,e\,o\,r\,e\,m \, 5. Let $M$ be an almost Hermitian manifold of indefinite metric
and $ n\ge 3$. If $R(X,\xi,\xi,X)=0$, whenever span$\{X,\xi\}$ is a weakly isotropic
antiholomorphic 2-plane, then $M$ is of pointwise constant antiholomorphic sectional
curvature.

{\it Proof.} Let $n-s\ge 2$ (the case $s\ge 2$ is treated analogously). We choose $x,\,y,\,a$
in $T_pM$, $p\in M$ so that $\{x,\,y,\,a\}$ is an orthonormal antiholomorphic triple of 
signature $(+,+,-)$. From the condition of the theorem we have
$$
	R(x,y+a,y+a,x)=0 \ ,
$$
From here we find
$$
	K(x,y)=K(x,a) \ ,   \leqno (4)
$$
$$
	R(x,y,a,x)=0   \ .   \leqno (5)
$$
Replacing $a$ by $(a+Ja)/\sqrt 2$ in (4) we get
$$
	R(x,a,Ja,x)=0   \ .   \leqno (6)
$$
Analogously
$$
	R(x,y,Jy,x)=0   \ .   \leqno (7)
$$
Now, let $Y,\,Z$ be arbitrary unit vectors in $T_pM$ so that $Y,\,Z\perp x,\ Jx$. If
$Y,\,Z$ have the same signature, then (4) implies
$$
	K(x,Y)=K(x,Z)   \leqno (8)
$$
in the following way. If $\{Y,\,Z\}$ has a signature $(+,+)$, we choose 
$a\perp x,\,Jx,\,\,Y,\,JY,\,Z,$ $JZ$ and then $K(x,Y)=K(x,a)=K(x,Z)$ 
according to (4). If $\{Y,\,Z\}$ has a signature $(-,-)$, we choose 
$y\perp x,\,Jx,\,\,Y,\,JY,\,Z,\,JZ$ and further $K(x,Y)=K(x,y)=K(x,Z)$. If
$Y,\,Z$ are of different signature, then (4) - (7) imply (8), Similarly,
if $g(a,a)=-1$, we find 
$$
	K(a,Y)=K(a,Z)    \leqno (9)
$$
for unit vectors $Y,\,Z$ in $T_pM$ and $Y,\,Z \perp a,\,Ja$.

Finally, let $\alpha,\,\beta$ be arbitrary nondegenerate antiholomorphic 2-planes
with orthonormal bases $\{X,\,Y\}$ and $\{Z,\,U\}$, respectively, and let 
$E={\rm span}\{X,\,JX\}$. If $\{Z,U\}$ is of signature $(+,-)$, we choose a
vector $W \in T_pM,W\perp\ Z,\,JZ$, $g(W,W)=1$. Then $\{Z,\ W\}$ is of
signature (+,+) and 
$$
	K(\beta)=K(Z,U)=K(Z,W)
$$
taking into account (8). Hence we can assume $\{Z,\,U\}$ is of signature (+,+)
or $(-,-)$. Let $Z'$ be a unit vector and $Z' \in \beta \cap E^\perp$. 
Choosing $U' \in \beta$, $U'\perp Z'$ and using (8), (9), we obtain
$$
	K(\alpha)=K(X,Y)=K(X,Z')=K(Z',U')=K(\beta)  \ .
$$
This gives: $M$ is of pointwise constant antiholomorphic sectional curvature.

\vspace{0.1cm}
R\,e\,m\,a\,r\,k. \ If $M$ ($n\ge 3$) is an almost Hermitian manifold of pointwise 
constant antiholomorphic sectional curvature $\nu$, then its curvature tensor has the form [5]
$$
	R-\frac{1}{2(n+1)} \psi(\rho^*) + \frac{\tau^*}{(2n+1)(2n+2)}\pi_2
	=\nu\left(\pi_1- \frac{1}{2n+1}\pi_2\right) \ ,
$$
where
$$
	\begin{array}{r}
		\psi(S)(x,y,z,u)=g(y,Jz)S(x,Ju)-g(x,Jz)S(y,Ju)-2g(x,Jy)S(z,Ju) \\
		+g(x,Ju)S(y,Jz)-g(y,Ju)S(x,Jz)-2g(z,Ju)S(x,Jy)
	\end{array}
$$
is a curvature like-tensor whenever $S(x,Jy)+S(y,Jx)=0$ and 
$$
	\pi_2(x,y,z,u)=g(y,Jz)g(x,Ju)-g(x,Jz)g(y,Ju)-2g(x,Jy)g(z,Ju) \ .
$$
Then, the direct verification shows the inverse proposition of Theorem 5.

If $M$ is Kaehlerian and $n\ge 3$, the following conditions are equivalent [8]

1) $M$ is of constant holomorphic sectional curvature $\mu$;

2) $M$ is of constant antiholomorphic sectional curvature $\mu/4$.

Then Theorem 5 implies

\vspace{0.1cm}
C\,o\,r\,o\,l\,l\,a\,r\,y \ [6]. Let $M$ be a Kaehler manifold of indefinite 
metric and $n\ge 3$. If $R(X,\xi,\xi,X)=0$, whenever span$\{X,\xi\}$ is a weakly 
isotropic antiholomorphic 2-plane, then $M$ is of constant holomorphic sectional
curvature.

\vspace{0.2cm}
{\bf Almost Hermitian manifolds with vanishing Bochner curvature tensor.}
Firstly, we shall prove

\vspace{0.1cm}
L\,e\,m\,m\,a \ 1. Let $T$ be a curvature-like tensor of type (0,4) over $T_pM$. If

1) $T(x,Jx,Jx,x)=0$, whenever $g(x,x)=1$,

2) $T(x,a,a,x)=0$, whenever span$\{x,\,a\}$ is an antiholomorphic 2-plane of signature $(+,-)$,

3) $T(x,Jx,b,x)=0$, whenever span$\{x,\,b\}$ is an antiholomorphic 2-plane of signature $(+,-)$, \\
then $T=0$.
 
{\it Proof.} We shall show that the conditions 1), 2) and 3) imply the corresponding
conditions of Theorem C. For instance, let 1) hold good. If $\alpha \in T_pM$ and
$g(a,a)=-1$, we choose $x\perp a,\,Ja$ and $g(x,x)=1$. Then 1) implies
$$
	T(x+ta,Jx+tJa,Jx+tJa,x+ta)=0
$$
for every $|t|<1$. From here we find $T(a,Ja,Ja,a)=0$ and hence $T(X,JX,JX,X)=0$ for
arbitrary non-isotropic vector $X$. Further, approximating any isotropic vector $\xi$
with non-isotropic vectors, we obtain $T(\xi,J\xi,J\xi,\xi)=0$ and the condition 1)
of Theorem C. The conditions 2) and 3) in Theorem C follow in a similar way.

The Bochner curvature tensor $B(R)$ for an almost Hermitian manifold $M$ $(2n \ge 6)$
is given by the equality [7]
$$
	B(R)=R-\{16(n+2)\}^{-1}(\varphi+\psi)(\rho+3\rho^*)(R+\overline R)
$$
$$
	-\{16(n-2)\}^{-1}(3\varphi-\psi)(\rho-\rho^*)(R+\overline R)-\{(4(n+1))^{-1}\psi(\rho^*)
	-(4(n-1))^{-1}\varphi(\rho)\}(R-\overline R)
$$
$$
	+\{16(n+1)(n+2)\}^{-1}(\tau+3\tau^*)(R)(\pi_1+\pi_2)
$$
$$
	+\{16(n-1)(n-2)\}^{-1}(\tau-\tau^*)(R)(3\pi_1-\pi_2)  \ ,
$$
where $\overline R(x,y,z,u)=R(Jx,Jy,Jz,Ju)$ for all $x,\,y,\,z,\,u \in T_pM$.

The next theorem gives a characterization of almost Hermitian manifolds of indefinite
metric with vanishing Bochner curvature tensor.

\vspace{0.07cm}
T\,h\,e\,o\,r\,e\,m \, 6. Let $M$ ($s\ge2,\ n-s\ge2$) be an almost Hermitian 
manifold of indefinite metric. The following conditions are equivalent:

1) $R(\xi,\eta,\eta,\xi)=0$, whenever span$\{\xi,\,\eta\}$ is a strongly isotropic 
antiholomorphic 2-plane,

2) $B(R)=0$.

{\it Proof.} Let $\xi$ be an arbitrary isotropic vector in $T_pM$. Choosing an isotropic
vector $\eta$, so that $\eta \perp\xi,\,J\xi$, from the condition we find
$$
	R(\xi,J\xi+\eta,J\xi+\eta,\xi)=0 \ .
$$
From here we obtain
$$
	R(\xi,J\xi,J\xi,\xi)=0 \ .
$$
This equality gives
$$
	R(x+a,Jx+Ja,Jx+Ja,x+a)= 0
$$
for an arbitrary orthonormal pair $\{x,\,a\}$ of signature $(+,-)$. Using the
last equality we get
$$
	K(x,Jx)+K(a,Ja)=K(x,Ja)+K(Jx,a)-2R(x,Jx,Ja,a)-2R(x,Ja,Jx,a)  .   \leqno (10)
$$
Let $y$ be in $T_pM$, so that $\{x,\,y,\,a\}$ is an orthonormal antiholomorphic triple 
of signature $(+,+,-)$. Replacing $a$ by $(a+ty)/\sqrt{1-t^2}$, $|t|<1$ in (10) we find
$$
	K(x,Jx)+K(y,Jy)=K(x,Jy)+K(Jx,y)+2R(x,Jx,Jy,y)+2R(x,Jy,Jx,y)  .   \leqno (11)
$$
If $\{ e_1,\hdots,e_{2n}\}$ is an orthonormal basis of $T_pM$, from (10) and (11) we obtain
$$
	\sum_{i=1}^{2n} K(e_i,Je_i) = \frac{\tau+3\tau^*}{2(n+1)} \ ,  \leqno (12)
$$
$$
	K(x,Jx)=\frac1{2(n+2)}\{\rho(x,x)+\rho(Jx,Jx)+6\rho^*(x,x)\}-\frac{\tau+3\tau^*}{4(n+1)(n+2)} \ . \leqno (13)
$$
Now, let  $\{x,\,y,\,a,\,b\}$ be an orthonormal antiholomorphic quadruple of signature
$(+,+,-,-)$. We have 
$$
	R(x+a,y+b,y+b,x+a)=0
$$
and further
$$
	R(x,y+b,y+b,x)+R(a,y+b,y+b,a)=0 \ .
$$
From here, for an arbitrary $z\in T_pM$, $z\perp y,\,Jy,\,a,\,Ja,\,b,\,Jb$ we have
$$
	R(x,y+b,y+b,x)-R(z,y+b,y+b,z)=0 \ .
$$
The last two equalities imply
$$
	\begin{array}{c}
		2(n-2)R(x,y+b,y+b,x)=\rho(y+b,y+b)-R(Jy,y+b,y+b,Jy) \\
		+R(Jb,y+b,y+b,Jb) \ .
	\end{array}   \leqno (14)
$$
Firstly, from (14) it follows
$$
	\begin{array}{c}
		2(n-2)\{K(x,y)-K(x,b)\}=\rho(y,y)+\rho(b,b)-K(y,Jy)+K(b,Jb) \\
		+K(Jy,b)-K(y,Jb) \ .
	\end{array}   
$$
Using (12), from this equallity we obtain
$$
	(4n^2-14n+11)K(x,b)+(2n-3)K(x,Jb)+K(Jx,Jb)-K(Jx,b)
$$
$$
	+2(n-1)\{K(x,Jx)+K(b,Jb)\}=(2n-3)\rho(x,x)-2(n-2)\rho(b,b)
$$
$$
	+\rho(Jx,Jx)-2\rho(Jb,Jb)-\tau+\frac{\tau+3\tau^*}{2(n+1)} \ .
$$
This formula gives 
$$
	4(n-1)(n-2)K(x,b)+2(n-1)\{K(x,Jx)+K(b,Jb)\}
$$
$$
	=(2n-3)\{\rho(x,x)-\rho(b,b)\} +\rho(Jx,Jx)+\rho(Jb,Jb)-\tau+\frac{\tau+3\tau^*}{2(n+1)} \ .
$$
From here, taking into accout (13), we find
$$
	\begin{array}{c}
		\ds K(x,b)=\frac{2n^2-5}{4(n-1)(n^2-4)}\{\rho(x,x)-\rho(b,b)\} +\frac3{4(n-1)(n^2-4)}\{\rho(Jx,Jx) \\
		\ds -\rho(Jb,Jb)\}-\frac3{2(n^2-4)}\{\rho(x,x)-\rho(b,b)\}-\frac{2n^2+3n+4}{8(n^2-1)(n^2-4)}\tau  \\
		\ds + \frac{9n}{8(n^2-1)(n^2-4)}\tau^* \ .
	\end{array}   \leqno (15)
$$
Further, (14) implies
$$
	2(n-2)R(x,y,b,x)=\rho(y,b)-R(Jy,y,b,Jy)+R(Jb,y,b,Jb) \ .
$$
Replacing here $x$ by $(Jx-Jy)/\sqrt 2$ and $y$ by $(x+y)/\sqrt2$ we get
$$
	2(n-2)R(y,Jy,Jy,b)-(2n-5)\{R(x,Jx,Jy,b)+R(Jy,x,b,Jx)\}
$$
$$
	-R(b,Jb,Jb,y)+R(Jx,y,b,Jx)=\rho(y,b) \ ,
$$
from where it follows
$$
	\begin{array}{c}
		2nR(y,Jy,Jy,b)-R(b,Jb,Jb,y)-3R(b,Jb,y,b) \\
		=\rho(y,b)+3\rho^*(b,y) \ .    
	\end{array}   \leqno (16)
$$
for an arbitrary orthonormal pair $\{y,b\}$ of signature $(+,-)$.

Changing in (15) \ $x$ by $(y+tb)/\sqrt{1-t^2}$ and \ $b$ \ by $(tJy+Jb)/\sqrt{1-t^2}$
for arbitrary $|t|<1$ we check
$$
2(n-1)\{R(Jy,y,y,Jb)+R(Jb,y,b,Jb)\} = -\rho(y,b)+\rho(Jy,Jb) \ .  \leqno (17)
$$ 
Finally, from
$$
	R(y+b,Jy+Jb,Jy+Jb,y+b)=0
$$
we derive
$$
	R(y,Jy,Jy,b)+R(b,Jb,Jb,y)+R(y,Jy,Jb,y)+R(b,Jb,Jy,b)=0 \ .   \leqno (18)
$$
From (16), (17) and (18) it follows easily
$$
	\begin{array}{c}
		\ds R(y,Jy,Jy,b)=-\frac{3}{4(n-1)(n+2)}\rho(y,b) + \frac{2n+1}{4(n-1)(n+2)}\rho(Jy,Jb) \\
		\ds - \frac3{4(n+1)(n+2)}\rho^*(b,y) + \frac{3(2n+3)}{4(n+1)(n+2)}\rho^*(y,b) \ .    
	\end{array}   \leqno (19)
$$
Taking into account (13), (15) and (19) and Lemma 1 we obtain $B(R)=0$.

The inverse is a straightforward verification.

\vspace{0.07cm}
L\,e\,m\,m\,a \ 2. Let $M$ ($s\ge 2,\ n-s\ge2$) be a Kaehlerian manifold of 
indefinite metric. The following conditions are equivalent:

1) $R(\xi,J\xi,J\xi,\xi)=0$, whenever $\xi$ is an isotropic vector,

2) $R(\xi,\eta,\eta,\eta,\xi)=0$, whenever span$\{\xi,\,\eta\}$ is a strongly isotropic
antiholomorphic 2-plane. 

{\it Proof.} The implication 2) $\rightarrow$ 1) was shown in the proof of Theorem 6 
in the case  of an arbitrary almost Hermitian manifold.

Now, let $M$ satisfy 1).  This condition gives
$$
	R(\xi+\eta, J\xi+J\eta,J\xi+J\eta,\xi+\eta)=0 \ ,
$$
which implies
$$
	2R(\xi,J\xi,J\eta,\eta)+2R(\xi,J\eta,J\xi,\eta)+R(\xi,J\eta,J\eta,\xi)+R(\eta,J\xi,J\xi,\eta)=0 \ .
$$
Using this equality and the first Bianchi identity we find
$$
	R(\xi,\eta,\eta,\xi)+3R(\xi,J\eta,J\eta,\xi)=0 \ .
$$
Now 2) follows in a straightforward way.

\vspace{0.07cm}
Applying Lemma 2 and Theorem 6 we obtain

\vspace{0.07cm}
T\,h\,e\,o\,r\,e\,m \, 7. Let $M$ ($s\ge 2,\ n-s\ge2$) be a Kaehlerian manifold of 
indefinite metric. The following conditions are equivalent:

1) $R(\xi,J\xi,J\xi,\xi)=0$, whenever $\xi$ is an isotropic vector,

2) $B(R)=0$.

\vspace{0.1cm}
R\,e\,m\,a\,r\,k. \, In [8] there is announced the following proposition without a proof:   

\vspace{0.1cm}
Theorem  6.5. Let $M$ be a connected indefinite Kaehler manifold with 
complex dimension $n\ge 2$ and index $2s>0$. If $R(u,Ju,Ju,u)=0$ for all isotropic
vectors $u\in T_pM$ holds, then $M$ has constant holomorphic sectional curvature.

\vspace{0.2cm}
{\bf Plane axioms for almost Hermitian manifolds of indefinite metrics.} Let $M$
be an almost Hermitian manifold of indefinite metric. $M$ is said to satisfy the axiom
of the weakly (strongly) isotropic antiholomorphic 2-planes, if, for each point $p$
in $M$ and for any weakly (strongly) isotropic antiholomorphic 2-plane $E$ in 
$T_pM$, there exists a totally geodesic submanifold $M'$ of $M$ containing $p$,
such that $T_pM'=E$ and all the tangent spaces of $M'$ are weakly (strongly) isotropic
antiholomorphic 2-planes in $M$. The axiom of isotropic holomorphic 2-planes is
formulated in a similar way. We note, that every isotropic holomorphic 2-plane
is necessarily strongly isotropic. 

In this section we discuss these axioms. We need the following propositions.

\vspace{0.07cm}
Theorem  D. Let $M$ be an almost Hermitian manifold of indefinite metric. If
$$
	R(x,Jx,Jx,a)=0  \ ,
$$
whenever $\{x,\,a\}$ is an orthonormal antiholomorphic pair of signature $(+,-)$, 
then $M$ is of pointwise constant holomorphic sectional curvature and $R=\overline R$.

\vspace{0.07cm}
This proposition can be derived in a similar way as the corresponding theorem in the
case of an almost Hermitian manifold of definite metric [9].

Using Lemma 1, we obtain immediately.

\vspace{0.07cm}L\,e\,m\,m\,a \ 3. Let $M$ ($s\ge 2,\, n-s\ge2$) be an almost Hermitian manifold of 
indefinite metric. If $M$ is of pointwise constant holomorphic sectional curvature
$\mu$, pointwise constant antiholomorphic sectional curvature $\nu$, and
$R=\overline R$, then the curvature tensor has the form:
$$
	R=\nu\pi_1+\frac{\mu-\nu}3\pi_2 \ .   \leqno (20)
$$

\vspace{0.07cm}Theorem  E. [7] Let $M$ be an almost Hermitian manifold and $2n\ge6$. If
the curvature tensor of $M$ has the form (20), then $M$ is of constant sectional 
curvature or $M$ is a Kaehler manifold of constant holomorphic sectional curvature.

\vspace{0.07cm}
The proof of this theorem is also valid in the case of an indefinite metric.

\vspace{0.07cm}
T\,h\,e\,o\,r\,e\,m \, 8. Let $M$ ($s\ge2,\ n-s\ge2$) be an almost Hermitian 
manifold of indefinite metric. If $M$ satisfies the axiom of weakly isotropic 
antiholomorphic 2-planes, then $M$ is of constant sectional curvature or $M$ 
is a Kaehler manifold of constant holomorphic sectional curvature.

{\it Proof.} Let $\{\xi, \,x\}$ be a basis of an weakly isotropic antiholomorphic 
2-plane in $T_pM$ with $g(\xi,\xi)=g(\xi,x)=g(\xi,Jx)=0$ and $g(x,x)=\pm 1$. By
the conditions of the theorem, there exists a totally geodesic submanifold $M'$
through $p$ such that $T_pM'={\rm span}\{\xi,\,x\}$. This implies $R(\xi,x,)x$ 
is in $T_pM'$. Hence $R(\xi,x,x,\xi)=0$. Applying Theorem 5, we obtain $M$ is of
pointwise constant antiholomorphic sectional curvature. On the other hand, since
$Jx\perp \xi,\,x$, then $R(\xi,x,x,Jx)=0$. From here it follows immediately 
$R(Jx,x,x,a)=0$, whenever $\{x,\,a\}$ is an orthonormal antiholomorphic pair
of signature $(+,-)$. Applying Theorem D, Lemma 3 and Theorem E we obtain the
proposition.

\vspace{0.07cm}
T\,h\,e\,o\,r\,e\,m \, 9. Let $M$ ($s\ge2,\ n-s\ge2$) be an almost Hermitian 
manifold of indefinite metric. If $M$ satisfies the axiom of strongly isotropic 
antiholomorphic 2-planes, then $M$ has vanishing Bochner curvature tensor.

{\it Proof.} Let $\{\xi, \,\eta\}$ be a basis of a strongly isotropic antiholomorphic
2-plane in $T_pM$. From the condition of the assertion we find $R(\xi,\eta)\eta$ is
in span$\{\xi,\,\eta\}$. Hence, $R(\xi,\eta,\eta,\xi)=0$. Applying Theorem 6
we obtain $B(R)=0$.

\vspace{0.1cm}
T\,h\,e\,o\,r\,e\,m \, 10. Let $M$ ($s\ge2,\ n-s\ge2$) be a Kaehlerian 
manifold of indefinite metric. If $M$ satisfies the axiom of  isotropic 
holomorphic 2-planes, then $M$ has vanishing Bochner curvature tensor.

{\it Proof.} Let $\xi$ be an arbitrary isotropic vector in $T_pM$. By the conditions
of the theorem we have $R(\xi,J\xi)J\xi$ is in span$\{\xi,\,J\xi\}$. Hence,
$R(\xi,J\xi,J\xi,\xi)=0$. Now, Theorem 7 implies the assertion.

\vspace{0.1cm}
T\,h\,e\,o\,r\,e\,m \, 11. Let $M$ ($s\ge2,\ n-s\ge2$) be a Kaehlerian 
manifold of indefinite metric. If $M$ satisfies the axiom of strongly isotropic 
antiholomorphic 2-planes, then $M$ is of constant holomorphic sectional curvature.

{\it Proof.} Applying Theorem 9 we find $B(R)=0$. Let $\{\xi, \,\eta\}$ be a basis 
of an arbitrary strongly isotropic antiholomorphic 2-plane in $T_pM$. The condition
of the theorem gives $R(\xi,\eta)\eta$ is a linear combination of $\xi$ and $\eta$.
From this condition and $B(R)=0$ we check $\rho(\xi,\xi)=0$ for an arbitrary
isotropic vector $\xi$. This condition implies $M$ is Einsteinian [10]. Hence $M$ is 
of constant holomorphic sectional curvature.

\vspace{0.7in}
\centerline{\large REFERENCES}

\vspace{0.1in}

\noindent
\ 1. D\,a\,j\,c\,z\,e\,r M., K. N\,o\,m\,i\,z\,u. On sectional curvature of indefinite metrics. -
 Math. 

\ \ \ Ann., {\bf 247}, 1980, 279-282.

\noindent
\ 2. S\,c\,h\,o\,u\,t\,e\,n, J. A. Ricci Calculus. Berlin - New York, Springer, 1954.
 
\noindent
\ 3. K\,u\,l\,k\,a\,r\,n\,i, R. S. Curvature structures and conformal transformations.
- J. Diff. 

\ \ \ Geom., {\bf 4}, 1970, 425-451.

\noindent 
\ 4. G\,a\,n\,c\,h\,e\,v, G. Almost Hermitian manifolds similar to the complex space forms. - C. 

\ \ \ R. Acad. Bulg. Sci., {\bf 32}, 1979, 1179-1182.

\noindent 
\ 5. G\,a\,n\,c\,h\,e\,v, G. On Bochner curvature tensor in almost Hermitian manifolds. - Pliska,

\ \ \ Studia mathematica bulgarica, {\bf 9}, 1987, 33-43.

\noindent
\ 6. C\,a\,r\,m\,e\,n, M., A. R\,o\,m\,e\,r\,o. Sur les vari\'et\'es k\"ahleriennes ind\'efinies 
\`a courbure 

\ \ \ sectionnelle holomorphe constante. - C. R. Acad. Paris, S\'er. A, {\bf 297},
1983, 343-344.

\noindent
\ 7. T\,r\,i\,c\,e\,r\,r\,i, F., L. Vanhecke. Curvature tensors on almost Hermitian manifolds. -

\ \ \ Trans. Amer. Math. Soc., {\bf 267}, 1981, 365-394. 

\noindent
\ 8. B\,a\,r\,r\,o\,s, M., A. R\,o\,m\,e\,r\,o. Indefinite Kaehlerian manifolds. 
Math. Ann., {\bf 261}, 1982, 

\ \ \ 55-62.

\noindent
\ 9. K\,a\,s\,s\,a\,b\,o\,v, O. On the axiom of planes and the axiom of spheres in the 
almost Hermi-

\ \ \ tian geometry. - Serdica, {\bf 8}, 1982, 109-114.

\noindent
10. D\,a\,j\,c\,z\,e\,r, M., K. N\,o\,m\,i\,z\,u. On the boundedness of Ricci curvature of an
indefinite 

\ \ \ metric. -  Bol. Soc. Brasil. Mat., {\bf 11}, 1980, 25-30.

\vspace {0.3in}

\ \ \ \ \ \ \ \ \ \ \ \ \ \ \
Received 10.04.1985

\end{document}